\documentclass[reqno,letterpaper]{amsart}

\usepackage{hyperref}
\usepackage{amssymb}
\usepackage{amsthm}

\newcounter{mnotecount}[section]

\title
[Complex potentials and trapping]
{A decay estimate for a wave equation with trapping and a complex
  potential}

\author[L. Andersson]{Lars Andersson}
\email{laan@aei.mpg.de}
\address{Albert Einstein Institute, Am M\"uhlenberg 1, D-14476 Potsdam,
  Germany}

\author[P. Blue]{Pieter Blue}
\email{P.Blue@ed.ac.uk}
\address{The School of Mathematics and the Maxwell Institute,  University of Edinburgh,James Clerk Maxwell Building, 
The King's Buildings, 
Mayfield Road, 
Edinburgh, 
Scotland EH9 3JZ, UK}

\author[J.-P. Nicolas]{Jean-Philippe Nicolas}
\email{Jean-Philippe.Nicolas@univ-brest.fr}
\address{Laboratoire de Math\'ematiques, Universit\'e de Brest, 
6 avenue Victor Le Gorgeu, 29200 Brest,
    France}

\date{\today \ {\em File:\jobname{.tex}}}

\newtheorem{theorem}{Theorem}
\newtheorem{remark}[theorem]{Remark}

\newcommand{\Reals}{\mathbb{R}}

\newcommand{\Sphere}{S}

\renewcommand{\Re}{\mathrm{Re}}
\renewcommand{\Im}{\mathrm{Im}}

\newcommand{\pt}{\partial_t}
\newcommand{\px}{\partial_x}
\newcommand{\pomega}{\nabla_\omega}
\newcommand{\Slap}{\Delta_{\omega}}

\newcommand{\di}{\mathrm{d}}
\newcommand{\dt}{\di t}
\newcommand{\dx}{\di x}
\newcommand{\domega}{\di^2\omega}
\newcommand{\dthree}{\,\,\dx\domega}
\newcommand{\dfour}{\,\,\dt\dx\domega}

\newcommand{\spacetime}{M}
\newcommand{\slab}[1]{M_{#1}}

\newcommand{\potlr}{V}
\newcommand{\potli}{W}
\newcommand{\wave}{u}
\newcommand{\wavebar}{\bar{u}}
\newcommand{\waver}{v}
\newcommand{\wavei}{w}

\newcommand{\initdata}{\psi_0}
\newcommand{\initdatat}{\psi_1}

\newcommand{\LagrangianA}{\mathcal{L}_2}
\newcommand{\LagrangianB}{\mathcal{L}_1}

\newcommand{\Energy}{E}
\newcommand{\GenEnergy}[1]{E_{#1}}

\begin{document}

\begin{abstract}
In this brief note, we consider a wave equation that has both trapping
and a complex potential. For this problem, we prove a uniform bound on
the energy and a Morawetz (or integrated local energy decay)
estimate. The equation is a model problem for certain scalar equations appearing
in the Maxwell and linearised Einstein systems on the exterior of a rotating
black hole.
\end{abstract}

\maketitle

\section{Introduction}
We consider the Cauchy problem 
\begin{align}\label{eq:Model} 
\left(-\pt^2+\px^2+\potlr(\Slap-N)   +i\epsilon\potli \right) \wave = 0 ,
\\
\wave(0,x,\omega) =\initdata(x,\omega), \quad 
\pt\wave(0,x,\omega) =\initdatat&(x,\omega),  
\end{align} 
on
$(t,x,\omega)\in\spacetime=\Reals\times\Reals\times\Sphere^2$ 
with smooth, compactly supported initial data. 
Here $\wave$ is a complex function 
$\wave=\waver+i\wavei$, 
\begin{align*}
\potlr&= \frac{1}{x^2+1} , 
\end{align*}
$\potli$ is a smooth, real-valued, compactly supported function
which is
nonvanishing at $x=0$ and uniformly bounded by $1$, and $\epsilon>0$
is a small parameter. Finally, $\Slap$ is the Laplacian in the
angular variables and $N$ is a number chosen to be sufficiently large
to allow us to avoid certain technical problems.

The equation \eqref{eq:Model} has both trapping, which occurs at
$x=0$, and a complex potential, which does not vanish at the trapped
set. The interaction of these creates problems, which appear to
frustrate the use of energy and Morawetz estimates at the classical
level. 
By adapting known, pseudodifferential methods, we
show how to overcome these problems. 
We now state our main result in terms of the energy
\begin{align*}
\Energy(t) &= \frac12\int_{\{t\}\times\Reals\times\Sphere^2} 
|\pt\wave|^2 +|\px\wave|^2 +\potlr(|\pomega\wave|^2+N|\wave|^2)
\dthree ,
\end{align*}

\begin{theorem}
\label{thm:Main}
There is a constant $C$ such that, if $\initdata$ and $\initdatat$ are such that $\Energy(0)$ is finite then 
\begin{subequations}
\begin{gather}
\forall t\in\Reals:\hspace{2cm}
\Energy(t)\leq C\Energy(0) , 
\label{eq:EnergyBound}\\
\int_{\spacetime}
\frac{|\px\wave|^2}{x^2+1}
+\frac{x^2}{1+x^2}\left(\frac{|\pomega\wave|^2}{1+|x|^3} +\frac{|\pt\wave|^2}{x^2+1}\right)
+\frac{|\wave|^2}{1+|x|^3} 
\dfour \leq C\Energy(0) , 
\label{eq:ClassicalMorawetz}\\
\int_{\spacetime}
\frac{|u| |\pt\wave|}{1+|x|^3} \dfour \leq C\Energy(0) . 
\label{eq:RefinedMorawetz}
\end{gather}
\end{subequations}
\end{theorem}

Since the equation \eqref{eq:Model} has $t$ independent coefficients,
one might naively think that Noether's theorem provides a positive conserved
energy. However, 
for the Lagrangian $\LagrangianB[\wave,\partial\wave]
= -(\pt\wave)^2+(\px\wave)^2
+\potlr(\pomega\wave\cdot\pomega\wave+N\wave^2) - i\epsilon\potli\wave^2$,
which has the wave equation \eqref{eq:Model} as its
Euler-Lagrange equation, the conserved quantity associated to the
time translation symmetry is indefinite, being approximately the
energy of the real component of $\wave$ minus the energy of the
imaginary component (plus $\epsilon$ times a term involving
$\potli\waver\wavei$). On the other hand, a Lagrangian of the form
$\LagrangianA[\wave,\partial\wave] = -|\pt\wave|^2+|\px\wave|^2
+\potlr(|\pomega\wave|^2+N|\wave|^2)$, which corresponds to the energy
expression considered above, does not yield equation \eqref{eq:Model} as its
Euler-Lagrange equation. 

The wave equation \eqref{eq:Model} is a model for equations arising in
the study of the Maxwell and linearised Einstein equations outside a
Kerr black hole. The Kerr black holes are a family of Lorentzian
manifolds arising in general relativity, and they are characterised by
a mass parameter $M$ and an angular momentum parameter $a$. Black
holes are believed to be the enormously massive objects at the center
of most galaxies. The case $|a|\leq M$ is the physically relevant. The
case $a=0$ is the Schwarzschild class of black holes. 

It is expected that every uncharged black hole will asymptotically
approach a Kerr solution under the dynamics generated by the Einstein
equations of general relativity. The wave, Maxwell, 
and linearised
Einstein equations on a fixed Kerr geometry are a sequence of
increasingly accurate models for these dynamics. By projecting on a null
tetrad, the Maxwell and
linearised Einstein fields can be decomposed into sets of complex scalars,
the Newman-Penrose (NP) scalars
\cite{newman:penrose:1962,newman:penrose:errata,GHP}.  
It is well-known that the NP scalars with extreme spin weights 
satisfy
decoupled wave equations, known as the Teukolsky equations, and that the
solutions to these reduced equations can be used to reconstruct the
full system \cite{Teukolsky}. 

For the Maxwell field on the Kerr background,
the spin weight 0 NP scalar can be treated in the same way, and the resulting
equation is known as the Fackerell-Ipser equation \cite{FackerellIpser}. For 
linearized gravity on the Schwarzschild
background, it is also well known that the imaginary part of the 
spin weight 0 NP scalar
is governed by a wave equation, the 
Regge-Wheeler equation \cite{ReggeWheeler,price:1972:II}. The corresponding
equation for the real part is more complicated,
cf. \cite{zerilli:1970,moncrief:1974}, see also \cite{AksteinerAndersson}.

It was recently shown \cite{AksteinerAndersson} that
in the general $(|a|<M)$ Kerr case, 
by imposing a gauge
condition 
related to the wave coordinates gauge, the equation for both the real and
imaginary parts of the spin
weight 0 NP scalar of the linearized gravitational field may be put in a form
analogous to the Regge-Wheeler and Fackerell-Ipser equations. 
Explicitely (in the Kerr spacetime with signature $-$+++, working in Boyer-Lindquist
coordinates) these take the form 
\begin{equation}
\left( \nabla^\alpha \nabla_\alpha + 2s^2 \frac{M}{(r-ia\cos\theta)^3} \right) u = 0,
\label{eq:RW} 
\end{equation} 
where $s=0$ corresponds to the free scalar wave equation, 
$s=1$ corresponds to the Maxwell (Fackerell-Ipser) case, while $s=2$
corresponds to the linearized gravity (generalized Regge-Wheeler) case. 
In particular, for the $a\not=0$ cases, the
analogues of the Regge-Wheeler equations have complex potentials, with
the imaginary part depending continuously on $a$.

For the wave equation in the Schwarzschild case, the use of the energy
estimate \cite{Wald}, like \eqref{eq:EnergyBound} with $C=1$, and
Morawetz estimates\footnote{These are also called integrated local energy decay
estimates.} are well established \cite{LabaSoffer, BlueSoffer,
BlueSoffer:Errata, BlueSterbenz,
DafermosRodnianski:RedShiftSchwarzschild}. 
In the Morawetz estimate
\eqref{eq:ClassicalMorawetz}, there is a loss of control of time and
angular derivatives near $x=0$, in the sense that the integrand cannot
control
$|x|^{p}(|\wave|^{q_t}|\pt\wave|^{2-q_t}+|\wave|^{q_\omega}|\pomega\wave|^{2-q_\omega})$
with both $p=0$ and either $q_t=0$ or $q_\omega=0$. The presence of
trapping makes some loss unavoidable \cite{Ralston}. By applying
``angular modulation'' and ``phase space analysis'', the range for the
angular parameter $q_\omega$ can be refined to $p=0$ and $q>0$
\cite{BlueSoffer:LongPaper}. This type 
of refinement is crucial in the
current paper. Alternatively, certain pseudodifferential operators
have been used to obtain refinements near $x=0$, to $p>0$ and
$q_t=q_\omega=0$ \cite{MarzuolaMetcalfeTataruTohaneanu}.

For the wave equation in the general ($|a|<M$) Kerr case, it is
possible to apply Fourier transforms first in the $\phi$ and $t$
variables\footnote{Here $\phi$ is the azimuthal angle, which would be
  one component of $\omega$ in the notation of this paper.} and then
the remaining variables. The individual $\phi$ modes decay pointwise
\cite{FinsterKamranSmollerYau:KerrWaveDecay}. Although the problem has
a time-translation symmetry, because the generator of time
translations fails to be a time-like vector with respect to the
Lorentzian inner product of the Kerr geometry, there is no positive,
conserved energy. A major advance was the proof that, in the slowly
rotating case $|a|\ll M$, there is a uniform energy bound, like
estimate \eqref{eq:EnergyBound}. The first proof used an estimate
similar to \eqref{eq:ClassicalMorawetz}, but with additional
restrictions on the support of the Fourier transform
\cite{DafermosRodnianski:KerrEnergyBound}. Independent work
\cite{TataruTohaneanu} established estimates similar to
\eqref{eq:EnergyBound} and \eqref{eq:ClassicalMorawetz}, but with no
restriction on the Fourier support, and there were subsequent
pseudodifferential refinements \cite{Tohaneanu}. Also, the first two
authors have proved similar results using methods which require two
additional levels of regularity but which completely avoid the use of
Fourier transforms.  Morawetz estimates and refinements are a crucial
step in proving pointwise decay estimates
\cite{BlueSoffer:LongPaper,BlueSterbenz,DafermosRodnianski:RedShiftSchwarzschild,DafermosRodnianski:BHStability,Luk,Luk:Kerr}
and Strichartz estimates
\cite{MarzuolaMetcalfeTataruTohaneanu,Tohaneanu}, including the
long-conjectured, inverse-cubic, Price law
\cite{MetcalfeTataruTohaneanu:Price,Tataru:Price}.

The study of the Maxwell and linearised Einstein systems in the Kerr
geometry is still in its infancy. For the general Kerr case,
a certain transformed, separated version of the Teukolsky system
has no exponentially growing modes \cite{Whiting}. 
In the Schwarzschild case, the $\phi$ modes of the
Teukolsky equation decay pointwise
\cite{FinsterSmoller:HigherSpin}. 
Recently, improved decay
estimates for the Regge-Wheeler type equation \eqref{eq:RW} on the
Schwarzschild background, giving decay rates of $t^{-3}$, $t^{-4}$ and $t^{-6}$,
respectively, for $s=0,1,2$, 
have been proved \cite{donninger:etal:2009arXiv0911.3179D}. 
The Regge-Wheeler equation has been
used with the full system to prove energy and Morawetz (and pointwise
decay) estimates for the Maxwell system \cite{Blue:Maxwell} and, more
recently, under assumptions on the asymptotic behaviour, for the full
(not merely linearised) Einstein equation
\cite{Holzegel:SchwarzschildEinstein}.

For the spin-weight 0 equations arising from the Maxwell and
linearised Einstein equations with $a\not=0$, the presence of complex
potentials in the reduced equations \eqref{eq:RW} prevents the
existence of a positive conserved energy and means that an unrefined
Morawetz estimate, such as \eqref{eq:ClassicalMorawetz}, is
insufficient to control the growth of the energy, which is why we also
prove estimate \eqref{eq:RefinedMorawetz}. In the Kerr geometry, there
is no positive, conserved energy because the generator of the
time-translation symmetry fails to be timelike everywhere. In
contrast, for the equation \eqref{eq:Model}, there is no positive,
conserved energy because the complex potential prevents any
variational approach from providing an energy-momentum tensor that
satisfies the dominant energy condition. The method of this paper
combines a Fourier-transform-in-time technique (as in
\cite{DafermosRodnianski:KerrEnergyBound,TataruTohaneanu}) with a
``modulation'' (or Fourier-rescaling) technique (from
\cite{BlueSoffer:LongPaper}).

As is common, $C$ will be used to denote a constant which may vary
from line to line, but which is independent of the choice of $\wave$
or $T$. The notation $A\lesssim B$ is used to denote that there is
some $C$ such that $A< CB$, with $C$ independent of $\wave$ and $T$,
and similarly for $\gtrsim$.

\section{A preliminary energy estimate}
We derive an estimate for an energy for the wave equation
\eqref{eq:Model} by integrating by parts against $\pt\wavebar$ and
following the standard procedure for getting an energy estimate:
\begin{align*}
0
&=\Re\left( (\pt\wavebar) (-\pt^2+\px^2+\potlr(\Slap-N)+i\epsilon\potli)\wave \right)\\
&=-\frac12\pt |\pt\wave|^2 \\
&\quad +\px\Re((\pt\wavebar)\px\wave)-\frac12\pt|\px\wave|^2 \\
&\quad
+\pomega\cdot\Re((\pt\wavebar)\pomega\wave)-\frac12\pt\left(\potlr(|\pomega\wave|^2+N|\wave|^2)\right)
\\
&\quad -\epsilon\potli\Im((\pt\wavebar)\wave) . 
\end{align*}

Introducing an energy which we denote by 
\begin{align*}
\Energy(t) &= \frac12\int_{\{t\}\times\Reals\times\Sphere^2} 
|\pt\wave|^2 +|\px\wave|^2 +\potlr(|\pomega\wave|^2+N|\wave|^2)
\dthree ,
\end{align*}
assuming that $\wave$ decays sufficiently rapidly as
$|x|\rightarrow\infty$, and integrating the previous formula over a
region $[t_1,t_2]\times\Reals\times\Sphere^2$, we find
\begin{align}
E(t_2)-E(t_1)
&=\int_{[t_1,t_2]\times\Reals\times\Sphere^2}-\epsilon\potli\Im((\pt\wavebar)\wave)
\dfour .
\label{eq:EnergyEstimate}
\end{align}
In particular, note that the energy fails to be conserved and that an
estimate of the form \eqref{eq:ClassicalMorawetz} would be
insufficient to control the right-hand side. There is, however, a
trivial exponential bound:
\begin{align*}
E(t_2)\leq e^{\epsilon(t_2-t_1)}E(t_1) .
\end{align*}

\section{The Morawetz estimate}
Following the standard procedure for investigating the wave
equation, we derive a Morawetz
estimate by multiplying the wave equation by $(f(x)\px\wavebar +q(x)\wavebar)$, where $f$ and $q$ are real-valued functions. 

In performing this calculation, it is useful to observe that
\begin{align*}
q'(x)\Re(\wavebar\px\wave)
&= \px\left(\frac{q'}{2}\wavebar\wave\right)
-\frac12 q'' \wavebar\wave . 
\end{align*}
Using this and applying the product rule term-by-term, one finds
\begin{align}
\Re\big(\left(f\px\wavebar+q\wavebar\right)&\left(-\pt^2\wave+\px^2\wave+\potlr(\Slap -N)\wave+i\epsilon\potli\wave\right)\big)\nonumber\\
&=
\pt p_t
+\px p_x +\pomega p_\omega \nonumber\\
&\quad+\left(-\frac12f'+q\right)|\pt\wave|^2
-\left( \frac12f'+q\right)|\px\wave|^2\nonumber\\
&\quad+\left(\left(\frac12f'-q\right)\potlr+\frac12 f(\px\potlr)\right)|\pomega\wave|^2\nonumber\\
&\quad+\left(N\left(\left(\frac12f'-q\right)\potlr+\frac12 f(\px\potlr)\right)+\frac12q''\right)|\wave|^2 \nonumber\\
&\quad- \epsilon f\potli\Im((\px\wavebar)\wave) ,
\label{eq:MorawetzInDivergenceForm}
\end{align}
where 
\begin{align*}
p_t
&=p_t(f,q;\wave)=-\Re((f(\px\wavebar)+q\wavebar)(\pt\wave)), \\
p_x
&=p_x(f,q,\wave)=\frac12f|\pt\wave|^2+\frac12f|\px\wave|^2-\frac12 f\potlr|\pomega\wave|^2
\\
& \quad + q \Re ( \wavebar \px \wave )-\frac12(Nf\potlr+q')|\wave|^2 , \\
p_\omega
&=p_\omega(f,q;\wave)=f\potlr\Re((\px\wavebar)(\pomega\wave))+q\potlr\Re(\wavebar\pomega\wave) . 
\end{align*}

We take $f=-\arctan(x)$, for which $f'=-(x^2+1)^{-1}=-\potlr$,
$f''=2x(x^2+1)^{-2}$, and $f'''=-2(3x^2-1)(x^2+1)^{-3}$. We take
$q=f'/2+\delta(1+x^2)^{-1}\arctan(x)^2$ for some sufficiently small
$\delta$. 

We use the notation
\begin{align*}
\GenEnergy{f\px+q}(t)
&=\int_{\{t\}\times\Reals\times\Sphere^2} 
\Re( f (\px\wavebar)\pt\wave)+\Re(q\wavebar(\pt\wave))
\dthree ,
\end{align*}
and observe that, by a simple Cauchy-Schwarz argument, there is the
estimate $|\GenEnergy{f\px+q}|\leq C\Energy$.

Observing that the left-hand side of \eqref{eq:MorawetzInDivergenceForm} vanishes, we have
\begin{align*}
0
&=\pt p_t +\px p_x +\pomega p_\omega\\
&\quad+\delta\frac{\arctan(x)^2}{1+x^2}|\pt\wave|^2
+\frac{1}{1+x^2}(1-\delta\arctan(x)^2)|\px\wave|^2\\
&\quad+\left(\frac{x\arctan(x)-\delta\arctan(x)^2}{(1+x^2)^2}\right)|\pomega\wave|^2\nonumber \\
&\quad+\left(N\left(\frac{x\arctan(x)-\delta\arctan(x)^2}{(1+x^2)^2}\right)+\frac12q''\right)|\wave|^2\nonumber \\
&\quad- \epsilon f\potli\Im((\px\wavebar)\wave)
\end{align*}
Taking $\epsilon$ sufficiently small, $N$ sufficiently large, and
$\delta$ sufficiently small, the factors in front of $\vert \px \wave
\vert^2$ and $\vert \wave \vert^2$ are nonnegative and one can
dominate the term involving $W$ using these two terms. (These
estimates are uniform, in the sense that, if the estimate holds for
choices of $\epsilon_0$, $N_0$, and $\delta_0$, then it remains valid
for $\epsilon<\epsilon_0$, $N=N_0$, and $\delta=\delta_0$.) Thus, by
integrating over a time-space slab
$\slab{[t_1,t_2]}=[t_1,t_2]\times\Reals\times S^2$, one can conclude
that there is a constant $C$ such that
\begin{align}
&\Energy(t_2)+\Energy(t_1)\nonumber\\
&\gtrsim
\int_{\slab{[t_1,t_2]}}
\frac{|\px\wave|^2}{x^2+1}
+|\arctan(x)|^2\left(\frac{|\pomega\wave|^2}{1+|x|^3} +\frac{|\pt\wave|^2}{x^2+1}\right)
+\frac{|\wave|^2}{1+|x|^3} 
\dfour .
\label{eq:MorawetzStandard}
\end{align}

\section{Pseudodifferential refinements}

\subsection{The wave equation for an approximate solution}
\newcommand{\scale}{\alpha}
\newcommand{\cutofftA}{\chi_{1}}
\newcommand{\cutofftB}{\chi_{2}}
\newcommand{\cutofftC}{\chi_{3}}
\newcommand{\cutoffx}{\chi_{|x|\leq2}}

We define a smooth characteristic function of an interval $[a,b]$ to be a function which is identically $1$ on $[a,b]$, which is supported on $[a-1,b+1]$, and which is monotonic on each of the intervals $[a-1,a]$ and $[b,b+1]$. A smooth characteristic function of a collection of intervals, each of which are separated by distance at least two, is defined to be the sum of the smooth characteristic functions of each interval. 

Let $T>0$ be a large constant. (Here large means larger than $-\log |\epsilon|$ and $2$.) Let $\cutofftA$ be a smooth characteristic function on $[0,T]$, and let $\cutofftB$ be a smooth characteristic function of $[-1,0]\cup[T,T+1]$. Let $\cutoffx$ be a smooth characteristic function of $[-1,1]$. We will use $\cutofftA$, $\cutofftB$, and $\cutoffx$ to denote $\cutofftA(t)$, $\cutofftB(t)$, and $\cutoffx(x)$ respectively. 

Since $\cutofftA$ is smooth, there is a uniform bound on its derivative and second derivative, each of which are supported on $[0,1]\cup[T,T+1]$, so that there is a constant $C$ such that $|\partial_t\cutofftA|+|\partial_t^2\cutofftA|\leq C \cutofftB$. 

The functions
\newcommand{\approxsolu}{{u_1}}
\newcommand{\bdrysolu}{{u_2}}
\newcommand{\approxIntOnlysolu}{{u_3}}
\newcommand{\F}{F}
\newcommand{\G}{G}
\begin{align*}
\approxsolu          &=\cutofftA\cutoffx \wave ,\\
\bdrysolu            &=\cutofftB\cutoffx\wave ,\\
\approxIntOnlysolu   &=\cutofftA\wave, 
\end{align*}
satisfies the equation 
\begin{align}
\left(-\partial_t^2+\partial_x^2+\potlr(\Slap-N)+i\epsilon\potli\right)\approxsolu
&=\F(\bdrysolu,\nabla\bdrysolu,t,x) +\G(\approxIntOnlysolu,\nabla\approxIntOnlysolu,t,x), 
\label{eq:ApproxDivergence}
\end{align}
where
\begin{align*}
\F(\bdrysolu,\nabla\bdrysolu,t,x)
&=-2(\partial_t\cutofftA)(\partial_t\bdrysolu)-(\partial_t^2\cutofftA)\bdrysolu\\
\G(\approxIntOnlysolu,\nabla\approxIntOnlysolu,t,x)
&=2(\partial_x\cutoffx)(\partial_x\approxIntOnlysolu)+(\partial_x^2\cutoffx)\approxIntOnlysolu . 
\end{align*}
Since all functions of $t$ in this equation are smooth and supported in $t\in[-2,T+2]$, they are Schwartz class in $t$, so we may take the Fourier transform in $t$ and remain in the Schwartz class. We will use $\text{ }\widehat{}\text{ }$ to denote the Fourier transform in $t$, and $\tau$ for the argument of such functions. We will typically use the word ``functions'' to describe $\wave$, $\approxsolu$, $\bdrysolu$, and $\approxIntOnlysolu$ and the words ``Fourier transforms'' to describe their Fourier transforms. We will use $L^2$ to denote $L^2(\di\omega\di x\di t)$ for functions and to denote $L^2(\di\omega\di x\di\tau)$ for Fourier transforms. We will use $\|\cdot\|$ for $\|\cdot\|_{L^2}$ unless otherwise specified.  

We introduce the following space-time integrals
\begin{align*}
I(T)&=\int_{-2}^{T+2}\int_{-2}^{2}\int_{S^2} x^2|\pt\wave|^2+|\px\wave|^2 +|\wave|^2 \di\omega\di x\di t, \\
J(T)&=\int_{\Reals\times\Reals\times S^2}  |\tau|^{6/5} |\widehat{\approxsolu}|^2 \di\omega\di x\di\tau . 
\end{align*}
The dependence of $J$ upon $T$ is through the smooth cut-off
$\cutofftA$ in $\approxsolu$. Typically, the argument $T$ will be
clear from context and will be omitted. From the Morawetz estimate
\eqref{eq:MorawetzStandard} and the exponential bound on the energy,
it follows that $I\lesssim E(T)+E(0)$.

We now aim to prove a Morawetz estimate using the Fourier
transform. We take
\begin{align*}
f&=-\arctan(|\tau|^\scale x), \\
q&=\frac{f'}{2}=\frac12 \frac{|\tau|^\scale}{1+|\tau|^{2\scale}x^2} ,
\end{align*}
with $\scale\in[0,1/2]$. We multiply the Fourier transform of equation \eqref{eq:ApproxDivergence} by $(f\px+q)\bar{\widehat{\approxsolu}}$, and integrate the real part over $\Reals\times\Reals\times S^2$. This integral is convergent because all the functions are compactly supported in time, so the Fourier transforms are Schwartz class. 

\subsection{Controlling the terms arising from the cut-off}
We consider first the integral arising from the right-hand side of \eqref{eq:ApproxDivergence}. This is
\begin{align*}
\int_{\Reals\times\Reals\times S^2} \Re \left( \left((f\px+q)\bar{\widehat{\approxsolu}}\right) \left(\widehat{F}+\widehat{G}\right) \right) \di\omega\di x\di\tau 
&\leq \|(f\px+q)\bar{\widehat{\approxsolu}}\| \|\widehat{F}+\widehat{G}\| .
\end{align*}
The terms on the right can be estimated by 
\begin{align*}
\|(f\px+q)\bar{\widehat{\approxsolu}}\|
&\leq \|f\px\widehat{\approxsolu}\| +\|q\widehat{\approxsolu}\| , \\
\|f\px\widehat{\approxsolu}\|&
\lesssim \|\px\widehat{\approxsolu}\| \lesssim I^{1/2}, \\
\|q\widehat{\approxsolu}\|
&\lesssim \| |\tau|^\scale\widehat{\approxsolu}\|
\lesssim \|\widehat{\approxsolu}\| +\| |\tau|^{1/2}\widehat{\approxsolu}\| \lesssim I^{1/2} +J^{1/2} ,
\end{align*}
and
\begin{align*}
\|\widehat{F}+\widehat{G}\|
&\leq \|\widehat{F}\|+\|\widehat{G}\| .
\end{align*}
Because $G$ is supported only for $t\in[-1,T+1]$ and $x\in[-2,2]$, we have
\begin{align*}
\|\widehat{G}\|&\lesssim I^{1/2}. 
\end{align*}
Similarly, because $F$ is supported only for $t\in[-1,0]\cup[T,T+1]$ and $x\in[-2,2]$, we have that at each instant in $t$, the function $F$ is bounded in $L^2(\di x\di\omega)$ by either $CE(0)^{1/2}$ or $CE(T)^{1/2}$. Since we are considering two intervals in $t$ of length $1$, we have
\begin{align*}
\|\widehat{F}\|
&=\|F\|\lesssim C(E(0)^{1/2}+E(T)^{1/2}) .
\end{align*}
Thus, the terms on the right-hand side of the Fourier transform of \eqref{eq:ApproxDivergence} are bounded by
\begin{align}
\int_{\Reals\times\Reals\times S^2} \Re &\left((f\px+q)\bar{\widehat{\approxsolu}}\right) \left(\widehat{F}+\widehat{G}\right) \di\omega\di x\di\tau \nonumber\\
&\leq C(E(0)^{1/2}+E(T)^{1/2}+J^{1/2})(E(0)^{1/2}+E(T)^{1/2}) .
\label{eq:tauStrongerBdryTerms}
\end{align}

\subsection{The Morawetz estimate for the approximate solution}
If we multiply the left-hand side of the Fourier transform of the wave equation \eqref{eq:ApproxDivergence} by $(f\px+q)\bar{\widehat{\approxsolu}}$ and take the real part, then we have the analogue of \eqref{eq:MorawetzInDivergenceForm}
\begin{align}
\Re\big(\left(f\px\bar{\widehat{\approxsolu}}+q\bar{\widehat{\approxsolu}}\right)&\left(\tau^2\widehat{\approxsolu}+\px\widehat{\approxsolu}+\potlr(\Slap-N)\widehat{\approxsolu}+i\epsilon\potli\widehat{\approxsolu}\right)\big)\nonumber\\
&=
\px p_x +\pomega p_\omega \nonumber\\
&\quad+\left(-\frac12f'+q\right)|\tau\widehat{\approxsolu}|^2
-\left(\frac12f'+q\right)|\px\widehat{\approxsolu}|^2\nonumber\\
&\quad+\left(\left(\frac12f'-q\right)\potlr+\frac12 f(\px\potlr)\right)|\pomega\widehat{\approxsolu}|^2\nonumber\\
&\quad+\left(N\left(\left(\frac12f'-q\right)\potlr+\frac12 f(\px\potlr)\right)+\frac12q''\right)|\widehat{\approxsolu}|^2 \nonumber\\
&\quad- \epsilon f\potli\Im((\px\bar{\widehat{\approxsolu}})\widehat{\approxsolu}) ,
\end{align}
where 
\begin{align*}
p_x
&=p_x(f,q,\widehat{\approxsolu})=\frac12f|\pt\widehat{\approxsolu}|^2+\frac12f|\px\widehat{\approxsolu}|^2-\frac12f\potlr|\pomega\widehat{\approxsolu}|^2
\\
&\quad + q\Re(\bar{\widehat{\approxsolu}}\px\widehat\approxsolu)
-\frac12(Nf\potlr+q')|\widehat{\approxsolu}|^2 , \\
p_\omega
&=p_\omega(f,q;\widehat{\approxsolu})=fV\Re((\px\bar{\widehat{\approxsolu}})(\pomega\widehat{\approxsolu}))+q\potlr\Re(\bar{\widehat{\approxsolu}}\pomega\widehat{\approxsolu}) . 
\end{align*}
Note that there is no $p_t$ term because, for Fourier transforms, the analogue of the product rule is simply $-\overline{i\tau\widehat{\approxsolu}}\widehat{\approxsolu}=\bar{\widehat{\approxsolu}}i\tau\widehat{\approxsolu}$. 

When this equality is integrated over a space-time slab, the $p_x$ and
$p_\omega$ terms integrate to zero, and the remaining terms are all non-negative except for those arising from $q''$ and from $\potli$. The integral of the term involving $W$ is bounded by $I$. 

We now consider the term involving $q''$: 
\begin{align*}
\frac12q''|\widehat{\approxsolu}|^2
&=|\tau|^{3\scale}\frac{1-3|\tau|^{2\scale}x^2}{(1+|\tau|^{2\scale}x^2)^3}|\widehat{\approxsolu}|^2 .
\end{align*}
From the positivity of the remaining terms and the bound
\eqref{eq:tauStrongerBdryTerms} of the terms coming from the
right-hand side of the wave equation \eqref{eq:ApproxDivergence} for $\approxsolu$, we have that
\begin{align*}
\int_{\Reals\times\Reals\times S^2} &|\tau|^{3\scale}\frac{1-3|\tau|^{2\scale}x^2}{(1+|\tau|^{2\scale}x^2)^3}|\widehat{\approxsolu}|^2 \di\omega\di x\di\tau \\
&\leq C(E(0)^{1/2}+E(T)^{1/2}+J^{1/2})(E(0)^{1/2}+E(T)^{1/2}) .
\end{align*}
This can be combined with an additional factor of $MI$, where $M$ is a
large constant. ($700$ is sufficient.) The integral $I$ dominates the integral of $(|\tau|^2+1)x^2|\widehat{\approxsolu}|^2$ and is bounded by $C(E(T)+E(0))$. Thus, we have
\begin{align*}
\int_{\Reals\times\Reals\times S^2} &\left(|\tau|^{3\scale}\frac{1-3|\tau|^{2\scale}x^2}{(1+|\tau|^{2\scale}x^2)^3}+M(|\tau|^2+1)x^2\right)|\widehat{\approxsolu}|^2 \di\omega\di x\di\tau\\
\leq& C(E(0)^{1/2}+E(T)^{1/2}+J^{1/2})(E(0)^{1/2}+E(T)^{1/2}) .
\end{align*}
By considering the two cases $|\tau|^\scale|x|<M^{-1/2}$ and
$|\tau|^\scale|x|\geq M^{-1/2}$, one can see that if $2-2\alpha = 3\alpha$ (i.e. $\alpha=2/5$), then 
\begin{align*}
\left(|\tau|^{3\scale}\frac{1-3|\tau|^{2\scale}x^2}{(1+|\tau|^{2\scale}x^2)^3}+M(|\tau|^2+1)x^2\right)\cutoffx
&\geq C|\tau|^{3\scale}\cutoffx
\end{align*} 
and therefore we find
\begin{align*}
C(E(0)^{1/2}&+E(T)^{1/2}+J^{1/2})(E(0)^{1/2}+E(T)^{1/2}) \\
&\geq
\int_{\Reals\times\Reals\times S^2} \left(|\tau|^{3\scale}\frac{1-3|\tau|^{2\scale}x^2}{(1+|\tau|^{2\scale}x^2)^3}+M(|\tau|^2+1)x^2\right)|\widehat{\approxsolu}|^2 \di\omega\di x\di\tau \geq J, 
\end{align*}
which implies
\begin{align}
J\leq&C(E(T)+E(0)) .
\label{eq:StrongerMorawetz}
\end{align}

\subsection{Closing the energy estimate}
It is now possible to estimate the integral on the right-hand side of
the energy estimate \eqref{eq:EnergyEstimate}. For $|x|\geq 1$, the
right-hand side (using the compact support of $W$ and the
Cauchy-Schwarz estimate) can be dominated by $I \leq
C(E(T)+E(0))$. For $|x|\leq 1$, we would like to dominate the integral
over $\Reals\times\Reals \times S^2$ of
$|\potli\Im(\bar{\approxsolu}\pt\approxsolu)|$ by the integral
$J$. However, this is not entirely correct, because in $J$ there is a
contribution arising from the support of $\wave$ in the region
$t\in[-1,0]\cup[T,T+1]$. The error in this approximation is
bounded by $C(E(T)+E(0))$. Thus, we have
\begin{align*}
E(T)-E(0)
&\leq C\epsilon\left(E(T)+E(0)+\left|\int_{\Reals\times\Reals\times S^2} \potli\Im(\bar{\approxsolu}\pt\approxsolu) \di\omega\di x\di t\right|\right).
\end{align*}
We can also take the Fourier transform to obtain an estimate by
\begin{align*}
C\epsilon&\left(E(T)+E(0)+\left|\int_{\Reals\times\Reals\times S^2} \potli \, \Im ( \bar{\widehat{\approxsolu}}\widehat{\pt\approxsolu} ) \di\omega\di x\di t\right|\right)\\
&\leq C\epsilon\left(E(T)+E(0)+\int_{\Reals\times\Reals\times S^2} |\tau| \vert \widehat{\approxsolu}\vert^2 \di\omega\di x\di t\right) .
\end{align*}
The integrand is now controlled by $I+J$, which, from estimates \eqref{eq:MorawetzStandard} and \eqref{eq:StrongerMorawetz}, we know can be estimated also by the sum of the initial and final energies. This leaves the estimate
\begin{align*}
E(T)-E(0)
&\leq C\epsilon(E(T)+E(0)).
\end{align*}
By taking $\epsilon$ sufficiently small relative to the constant, we obtain a uniform bound on the energy
\begin{align*}
E(T) &\leq CE(0) .
\end{align*}
We note that, since all the constants were independent of $T$, the
estimate holds uniformly in $T$. This proves the first statement
\eqref{eq:EnergyBound}, in theorem \ref{thm:Main}. Combining this with
estimate \eqref{eq:MorawetzStandard} (and estimating
$x^2(1+x^2)\lesssim\arctan(x)^2$) gives the second,
\eqref{eq:ClassicalMorawetz}. Finally, the arguments of this section
and the bound on $I+J$, from estimates \eqref{eq:MorawetzStandard} and
\eqref{eq:StrongerMorawetz}, give the third result
\eqref{eq:RefinedMorawetz}.

\begin{remark}
Using the method given in \cite{BlueSoffer:LongPaper}, the stronger Morawetz estimate \eqref{eq:StrongerMorawetz} can be improved to control the integral of $|\tau|^{2-\varepsilon}|\widehat\wave|^2$ for any $\varepsilon>0$. Because of the presence of trapping, it is not possible to improve this to $|\tau|^2 |\widehat{\wave}|^2$. 
\end{remark}

\newcommand{\prd}{Phys. Rev. D} 
\bibliography{ModelProblem}
\bibliographystyle{abbrv}

\end{document}